\documentclass[amssymb,11pt]{amsart}
\usepackage{amscd,amssymb,euscript,amsthm}
\theoremstyle{plain}
\newtheorem{thm}{Theorem}[section] 
\newtheorem{cor}[thm]{Corollary}
\newtheorem{prop}[thm]{Proposition}
\newtheorem{lem}[thm]{Lemma}
\theoremstyle{definition}
\newtheorem{defn}[thm]{Definition}
\theoremstyle{remark}
\newtheorem{rem}[thm]{Remark}

\numberwithin{equation}{section}
\newcommand{\tr}{\operatorname{trace}}
\def\<{\left<}
\def\>{\right>}
\begin{document}
\title[Diagonals]{Diagonals of self-adjoint operators}
\author{William Arveson*}
\author{Richard V. Kadison}
\thanks{*supported by 
NSF grant DMS-0100487} 
\address{Department of Mathematics,
University of California, Berkeley, CA 94720\\
and Department of Mathematics, University 
of Pennsylvania, Philadelphia, PA 19104}
\email{arveson@mail.math.berkeley.edu, kadison@math.upenn.edu}
\date{15 October, 2002}

\begin{abstract}
The eigenvalues of a self-adjoint $n\times n$ matrix $A$ 
can be put into a decreasing sequence 
$\lambda=(\lambda_1,\dots,\lambda_n)$, with 
repetitions according to multiplicity, and  
the diagonal of $A$ 
is a point of ${\mathbb R}^n$ that bears some 
relation to $\lambda$.  The Schur-Horn theorem 
characterizes that relation in terms of a system 
of linear inequalities.

We prove an extension of the latter result 
for positive trace-class operators on infinite 
dimensional Hilbert spaces, generalizing 
results of one of us on the diagonals of 
projections.  
We also establish an appropriate counterpart of 
the Schur inequalities that relate spectral 
properties of self-adjoint operators 
in $II_1$ factors to their images under 
a conditional expectation onto a maximal 
abelian subalgebra.  
\end{abstract}
\maketitle

\section{Introduction}\label{S:in}

This paper presents some of the 
results of a project begun by the authors that 
is directed toward 
finding an appropriate common generalization of the Schur-Horn 
theorem (for matrices) to operators on an infinite-dimensional Hilbert 
space, and to operators in finite factors, in a form that 
would generalize 
work of one of us on projections in $II_1$ factors \cite{kadPnasI}, 
\cite{kadPnasII}.  

That 
project continues, and remains unfinished.  
The results below are satisfactory 
in the case of type $I$ factors, but are incomplete for 
finite factors.  
We are making these partial results public 
since there is renewed interest in these directions \cite{dSherman}, and 
it seems desirable to avoid duplication of effort.  Other aspects of this work 
were presented in Section 5 of \cite{kadSchur}.  

We 
point out that while the results of Section \ref{S:sec3} 
may appear to overlap 
with work of A. Neumann \cite{aNeumann1}, that is actually not the case.  
Neumann's results characterize the {\em closure} (in the $\ell^\infty$ norm) 
of the set of diagonals of self-adjoint operators with prescribed spectral properties. 
Here, on the other hand, we are concerned with the diagonals themselves, and 
not with their limits 
relative to any topology.  For example, one should compare 
Theorem 15 of \cite{kadPnasII} -- which characterizes the diagonals of projections --
with the corresponding result of \cite{aNeumann1} (Theorem 3.6 and Corollary 2.14)
to understand the extent to which 
subtlety is lost when one takes the closure relative to the $\ell^\infty$-norm.  
Along with the nature of the characterizations below, our methods  
also differ significantly from those of \cite{aNeumann1}.  
We thank Daniel Markiewicz 
for calling our attention to the paper \cite{aNeumann1} (also see \cite{aNeumann2}).  

\section{The Results of Schur and Horn}

We begin by stating 
the key assertion of Theorem 
5 from Alfred Horn's 1954 paper \cite{hornDiag}, 
which can be formulated as follows.   

\begin{thm}[Horn]\label{hornThm}
Let $\lambda_1\geq\cdots\geq\lambda_n$ 
and $p_1\geq\cdots\geq p_n$ be two decreasing 
sequences of real numbers 
satisfying $p_1+\cdots+p_k\leq \lambda_1+\cdots+\lambda_k$, 
$1\leq k\leq n-1$, and $p_1+\cdots+p_n=\lambda_1+\cdots+\lambda_n$.  
Then there is a self-adjoint $n\times n$ matrix with eigenvalues 
$\lambda_1,\dots,\lambda_n$ whose diagonal entries 
are $p_1,\dots,p_n$.  
\end{thm}

Another proof is offered at the end of this section.  
In a more coordinate-free formulation, Horn's theorem makes 
the following assertion.  Let $A$ be a self-adjoint operator 
on an $n$-dimensional Hilbert space $H$ with eigenvalues 
$\lambda_1\geq\cdots\geq \lambda_n$, and let $p_1\geq\cdots\geq p_n$ 
be a decreasing sequence that relates to $\lambda$ as in the 
hypothesis of Theorem \ref{hornThm}.  Then there is an 
orthonormal basis $e_1,\dots,e_n$ for $H$ such that 
$$
\langle Ae_k,e_k\rangle = p_k,\qquad k=1,\dots,n.
$$

The converse of Theorem \ref{hornThm} 
is also true, and this is the part of the 
composite Schur-Horn theorem that is attributed to Schur \cite{schurDiag}:
If there is a self-adjoint 
$n\times n$ matrix $A$ with eigenvalue 
sequence $\lambda=(\lambda_1,\dots,\lambda_n)$ 
with diagonal 
$p=(p_1,\dots,p_n)$, both written in decreasing order, 
then the inequalities 
\begin{equation}\label{Eq000}
p_1+\cdots+p_k\leq \lambda_1+\cdots+\lambda_k,\qquad 1\leq k\leq n
\end{equation}
of the hypothesis 
of Theorem \ref{hornThm} are satisfied, with 
equality holding for $k=n$.  That implication 
follows from classical estimates going back to 
Weyl \cite{weyl1}
(see the proof of Theorem \ref{thm1} below).  

There are other formulations of the Schur-Horn theorem that 
borrow from classical inequalities \cite{hlp}, 
the most notable one 
being the following.  Given a
sequence $\lambda=(\lambda_1,\dots,\lambda_n)$ of real 
numbers, let $\mathcal O_\lambda$ be the set of all 
$n\times n$ self-adjoint matrices having eigenvalue
sequence $\lambda$.  Then the set 
$E(\mathcal O_\lambda)$ of all diagonals 
of matrices in $\mathcal O_\lambda$ 
is the convex hull $\Lambda$ 
of the set of points $\lambda\circ\pi\in\mathbb R^n$, $\pi\in S_n$,  
obtained by permuting the components of $\lambda$.  
Schur's part of the Schur-Horn theorem becomes the assertion 
$E(\mathcal O_\lambda)\subseteq\Lambda$ while Theorem \ref{hornThm} implies  
$E(\mathcal O_\lambda)\supseteq \Lambda$.  
These formulations are discussed in \cite{hornDiag}.  

The Schur-Horn theorem has led to generalizations in 
several directions.  In 1973, Kostant \cite{KostantSHthm} put it into 
the context of actions of compact Lie groups 
(which generalize the unitary group $U(n)$).  Later Atiyah
\cite{atiyHorn}, 
and independently Guillemin and Sternberg 
\cite{GSI}, reformulated 
Kostant's result in 
the broader context of symplectic manifolds $M$ acted on by 
a torus $T$, and showed that for every moment map $\Phi$ for 
the $T$-action, the range of $\Phi$ is the convex hull of the 
images of the $T$-fixed points of $M$.  See \cite{knutExpo} 
for more detail.  

There is some connection between the finite dimensional 
Schur-Horn theorem \cite{hornDiag} 
and Horn's subsequent work on the eigenvalues of sums of matrices 
that culminated in the inequalities conjectured in 
\cite{hornEigSum}, as described 
in \cite{knutExpo}.  The Horn conjecture 
was recently proved, following work of 
Klyachko \cite{klyHorn} and the proof of the saturation 
conjecture by Knutson and Tao \cite{knutTao}.   
Thus, it may be appropriate to point out that Friedland 
has obtained a generalization 
of Klyachko's results to finite 
sums of positive trace-class 
operators acting on infinite 
dimensional Hilbert spaces \cite{friedl}.  

The purpose of this paper is to discuss two infinite dimensional 
formulations of the Schur-Horn Theorem.  
In Sections \ref{S:sec2}--\ref{S:sec3} we 
present a generalization of the Schur-Horn theorem to positive 
trace class operators on infinite dimensional Hilbert spaces.  
In Sections \ref{S:specDist}--\ref{S:sec6}, we reformulate these issues in 
the context of finite factors, and we establish appropriate 
versions of the Schur inequalities.  The $II_1$ version of 
Horn's result (Theorem \ref{hornThm}) is left 
as an open problem.

The first author wishes to thank Allen Knutson for helpful  
comments about the Schur-Horn theorem  
including the sketch of a ``calculus" proof, and 
for providing some key references.

\begin{proof}[Proof of Theorem \ref{hornThm}]
We show how one can deduce Theorem \ref{hornThm} 
from two results of \cite{kadPnasI}.  Let $p=(p_1,\dots,p_n)$ 
and $\lambda=(\lambda_1,\dots,\lambda_n)$ be two decreasing 
sequences satisfying the hypotheses of Theorem \ref{hornThm}.  
By lemma 5 of \cite{kadPnasI}, there is a sequence of points 
$x_1,\dots,x_n$ in $\mathbb R^n$ such that $x_1=\lambda$, 
$x_n=p$, and for $k=1,\dots, n-1$, $x_{k+1}$ can be 
expressed in terms of $x_k$ as follows
\begin{equation}\label{Eq0}
x_{k+1}=t_k\cdot x_k +  (1-t_k)\cdot x_k\circ\tau_k
\end{equation}
where $t_k$ is a number in the unit interval, $\tau_k$ is 
a transposition in $S_n$, and where 
$x\circ\tau$ denotes 
$(x_{\tau(1)},\dots,x_{\tau(n)})\in\mathbb R^n$.  

Given $x_1=\lambda,\dots,x_n=p\in\mathbb R^n$, 
$t_1,\dots,t_{n-1}\in[0,1]$ and 
$\tau_1,\dots,\tau_{n-1}\in S_n$ such that 
the relations (\ref{Eq0}) are satisfied,  we exhibit a 
sequence of self-adjoint matrices $A_1,\dots,A_n$ 
such that $A_k$ has eigenvalue list $\lambda$ and 
diagonal sequence $x_k$ as follows.  Theorem 6 of 
\cite{kadPnasI} asserts the following:  Given a 
self-adjoint $n\times n$ matrix $A=(a_{ij})$ 
with diagonal sequence 
$x$, and given a transposition $\tau$ in $S_n$ and a 
number $t\in[0,1]$, there is a unitary matrix 
$U$ such that the diagonal of $UAU^*$ is 
$t\cdot x+(1-t)\cdot x\circ\tau$.  The proof 
exhibits $U=(u_{ij})$ explicitly; 
if $\tau$ is the transposition $(ij)$ 
then $U$ coincides 
with the identity matrix except for the four terms 
$u_{ii}, u_{ij},u_{ji}, u_{jj}$ specified by
$$
\begin{pmatrix}
u_{ii}&u_{ij}\\
u_{ji}&u_{jj}
\end{pmatrix}
=
\begin{pmatrix}
z\cos\theta &\sin\theta& \\
-z\sin\theta&\cos\theta 
\end{pmatrix},
$$
where $z$ is a complex number of absolute value 
$1$ such that $za_{ij}$ is pure imaginary, and 
where $\theta$ satisfies $\cos^2\theta=t$.  
Let $A_1$ be the 
diagonal matrix with diagonal 
$\lambda=(\lambda_1,\dots,\lambda_n)=x_1$.  Given that 
$A_1, \cdots,A_k$ have been defined 
and satisfy the asserted conditions for $1\leq k<n$, 
the above result 
implies that there is 
a unitary matrix $U_k$ such that 
$$
\text{diag\,}(U_kA_kU_k^*)=t_{k}\cdot x_k + (1-t_k)\cdot x_k\circ\tau_k.  
$$ 
Setting $A_{k+1}=U_kA_kU_k^*$ and continuing inductively,   
we obtain a sequence of matrices 
$A_1,\dots,A_n$ whose last term 
$A_n=U_{n-1}\cdots U_1 A_1 U_1^*\cdots U_{n-1}^*$ is 
a self-adjoint 
matrix having eigenvalue list $\lambda$ and diagonal $p$.  
\end{proof}

\part{Type $I_\infty$ factors}

We first give a generalization of the Schur-Horn theorem 
to the case of positive trace-class operators acting 
on a separable infinite-dimensional Hilbert space.

\section{$\mathcal L^1$-closed unitary orbits}\label{S:sec2}

Let $H$ be a separable Hilbert space and let $A$ be a 
positive compact operator on $H$.  The sequence of 
eigenvalues of 
$A$ can be put into decreasing order, with repetitions 
according to the multiplicity of positive terms in 
the sequence, to obtain a sequence 
$\lambda=(\lambda_1,\lambda_2,\dots)$ 
satisfying $\lambda_1\geq\lambda_2\geq\cdots\geq0$, and 
we have 
$$
\lambda_1+\lambda_2+\cdots=\tr A \in[0,+\infty].  
$$
Such a decreasing sequence $\lambda$ will be 
called an {\em eigenvalue list}.  The preceding formula 
shows that $A$ is trace-class iff its eigenvalue list 
belongs to $\ell^1$, and the set 
of all eigenvalue lists in $\ell^1$ is a 
weak$^*$-closed cone, the weak$^*$-topology on $\ell^1$ 
arising from the identification of 
$\ell^1$ with the dual of $c_0$.  

The eigenvalue list of $A$ 
fails to be a complete invariant 
for unitary equivalence because it fails to 
detect zero eigenvalues except when $A$ is of finite rank.  For example, 
if $A$ has infinitely many positive 
terms $\lambda_k$ in its spectrum and has 
trivial kernel, then $A$ and $A\oplus \mathbf 0$ ($\mathbf 0$ being 
an the zero operator on some space of positive 
dimension) cannot 
be unitarily equivalent despite the fact that both have 
the same eigenvalue list.  

The state of affairs for trace-class operators 
is described as follows.  We write $\mathcal L^1=\mathcal L^1(H)$ for 
the Banach space of all trace-class operators on a Hilbert 
space $H$ with 
respect to the trace norm 
$$
\|A\|_1=\tr |A|, 
$$
$|A|$ denoting the positive square root of $A^*A$.  Given 
an eigenvalue list $\lambda\in\ell^1$, $\mathcal O_\lambda$ 
will denote the set of all positive trace-class operators on $H$ 
having $\lambda$ as their eigenvalue list.  Given a positive 
trace-class operator $A\in\mathcal B(H)$, $\mathcal O(A)$ 
will denote the {\em trace-norm closure} of the 
unitary orbit of $A$
$$
\mathcal O(A)=\{UAU^*: U\in\mathcal U(H)\}^{-\|\cdot\|_1}.  
$$
Two trace-class operators $A$, $B$ are said to be  $\mathcal L^1$- 
equivalent if there is a sequence of unitary operators $U_1, U_2,\dots$ 
such that 
that $\|A-U_nBU_n^*\|_1\to 0$ as $n\to\infty$; equivalently, 
$\mathcal O(A)=\mathcal O(B)$.  

\begin{prop}\label{prop1}
Let $A$ be a positive trace-class operator 
in $\mathcal B(H)$ and let $\lambda$ be the eigenvalue list of $A$.  
\begin{enumerate}
\item[(i)]
$\mathcal O(A)$ is a Polish topological space  
on which the unitary group of $H$ acts minimally.  
\item[(ii)]
$\mathcal O(A)=\mathcal O_\lambda$; in particular, the eigenvalue 
list is a complete invariant for $\mathcal L^1$-equivalence.  
\item[(iii)]
Two positive trace-class 
operators $A$, $B$ are $\mathcal L^1$-equivalent 
iff $A\oplus \mathbf 0$ and $B\oplus \mathbf 0$ are unitarily equivalent, 
where $\mathbf 0$ denotes the zero operator on an infinite 
dimensional Hilbert space.  
\item[(iv)]
If $\lambda$ has only finitely many nonzero terms, 
then $\mathcal O_\lambda$ consists of a single unitary 
orbit $\{UAU^*: U\in\mathcal U(H)\}$.  
\end{enumerate}
\end{prop}

\begin{proof}  (i): $\mathcal O_\lambda$ is a closed subset 
of $\mathcal L^1$ and therefore a separable complete metric 
space.  The fact that the orbit of every 
point of $\mathcal O(A)$ under the action 
of $\mathcal U(H)$ is dense 
in $\mathcal O(A)$ follows from the fact that 
$\mathcal L^1$-equivalence is a transitive relation.  

(ii):  Let $B$ be another positive trace-class operator 
with eigenvalue list $\mu$.  We have to show that 
$A$ and $B$ are $\mathcal L^1$-equivalent $\iff$ $\lambda=\mu$.  
For the implication $\implies$ we make use of the 
semiclassical inequality 
$$
\sum_{k=1}^\infty|\lambda_n-\mu_n|\leq \|A-B\|_1, 
$$
a proof of which can be found in the appendix of 
\cite{powThesis}.  Since $B$ can be closely approximated 
in the norm of $\mathcal L^1$ by operators unitarily 
equivalent to $A$, this inequality implies that 
$B$ must have the same eigenvalue list as $A$, hence   
$\mu=\lambda$.   Conversely, if $A$ and $B$ are 
two positive trace-class operators with the same 
eigenvalue list 
$\lambda = (\lambda_1,\lambda_2,\dots)$, then by the 
spectral theorem we can decompose $A$ and $B$ as follows
$$
A=A_n+R_n,\qquad B=B_n+S_n
$$
where $A_n$ and $B_n$ are finite rank positive operators 
with eigenvalue list $(\lambda_1,\dots,\lambda_n,0,0,\dots)$ 
and where the remainders $R_n$ and $S_n$ satisfy 
$$
\|R_n\|_1=\|S_n\|_1= \sum_{k=n+1}^\infty|\lambda_k|.  
$$
Since $A_n$ and $B_n$ are obviously 
unitarily equivalent for every 
$n=1,2,\dots$ and since $\|R_n\|_1$ and $\|S_n\|_1$ 
tend to zero as $n\to \infty$, it follows that 
there is a sequence of unitary operators 
$U_1, U_2,\dots$ such that 
$\|B-U_nAU_n^*\|_1\to 0$.  

(iii) is a consequence of (ii), which asserts  
that $A$ and $B$ are $\mathcal L^1$-equivalent iff 
they have the same eigenvalue list.  Indeed, it is obvious that 
if $A$ and $B$ have the same eigenvalue list $\lambda$ then 
$A\oplus\mathbf 0$ and $B\oplus\mathbf 0$ are unitarily 
equivalent; conversely, if $A\oplus\mathbf 0$ 
and $B\oplus\mathbf 0$ are unitarily equivalent then 
$A\oplus\mathbf 0$ and $B\oplus\mathbf 0$ 
must have the same eigenvalue list, hence so do $A$ and $B$.  

Finally, note that (iii)$\implies$(iv), since if $A$ is 
a finite rank positive operator with eigenvalue 
list $\lambda$, then all but a finite number of components of 
$\lambda$ are zero, hence $A$ is unitarily equivalent to 
$A\oplus\mathbf 0$, so that all operators in $\mathcal O_\lambda$ 
are unitarily equivalent.  
\end{proof}

\section{Diagonals of Trace Class Operators}\label{S:sec3}

Let $H$ be a separable 
Hilbert space and let 
$e_1, e_2,\dots$ be an orthonormal basis for $H$.  
The sequence of rank-one projections $E_k=[e_k]$, $k=1,2,\dots$ generates 
a discrete maximal abelian subalgebra 
$\mathcal A\subseteq\mathcal B(H)$, and the map that replaces 
an operator $A$ with the diagonal part 
$(a_{11}, a_{22},\dots)$ of its matrix 
$a_{ij}=\langle Ae_j,e_i\rangle$, 
$i,j=1,2,\dots$, relative to $(e_n)$ can be viewed as 
the unique trace preserving conditional 
expectation $E: \mathcal B(H)\to\mathcal A$ 
$$
E(A)=\sum_{n=1}^\infty E_nAE_n=\sum_{n=1}^\infty a_{nn} E_n.  
$$
The following result provides an infinite-dimensional 
generalization of the Schur-Horn theorem.  For a related result 
that characterizes the {\em norm-closure} of $E(\mathcal O(A))$
for a broader class of operators $A$, see \cite{aNeumann1}.  

\begin{thm}\label{thm1}Let $\mathcal A$ be a discrete 
maximal abelian von Neumann algebra in $\mathcal B(H)$, 
let $E:\mathcal B(H)\to\mathcal A$ be the 
trace-preserving conditional expectation on $\mathcal A$ 
and let $\lambda=(\lambda_1\geq\lambda_2\geq\cdots)$ be a decreasing 
sequence in $\ell^1$ with nonnegative terms.  
Then $E(\mathcal O_\lambda)$ consists of all positive trace-class 
operators $B\in\mathcal A$ whose eigenvalue list 
$p=(p_1\geq p_2\geq\cdots)$ satisfies 
\begin{equation}\label{Eq1}
p_1+\cdots+p_n\leq \lambda_1+\cdots+\lambda_n, \qquad n=1,2,\dots
\end{equation}
together with 
\begin{equation}\label{Eq2}
p_1+p_2+\cdots = \lambda_1+\lambda_2+\cdots.  
\end{equation}
\end{thm}

We will deduce Theorem \ref{thm1} from the following 
more general assertion about the diagonals of positive 
compact operators.  

\begin{thm}\label{thm2}
Let $\mathcal A\subseteq\mathcal B(H)$ be 
a discrete maximal abelian algebra with natural 
conditional expectation $E:\mathcal B(H)\to\mathcal A$.  
Let $A\in\mathcal B(H)$ be a 
positive compact operator with eigenvalue list 
$\lambda=(\lambda_1\geq \lambda_2\geq\cdots)$, 
and let $B$ be a positive compact operator in $\mathcal A$.  The 
following are equivalent.
\begin{enumerate}
\item[(i)] There is contraction $L\in\mathcal B(H)$ such 
that $E(L^*AL)=B$.  
\item[(ii)] The eigenvalue list $p=(p_1\geq p_2\geq\cdots)$ of 
$B$ satisfies 
$$
p_1+p_2+\cdots +p_n\leq \lambda_1+\lambda_2+\cdots+\lambda_n,\qquad 
n=1,2,\dots.  
$$
\end{enumerate}
\end{thm}

We require some eigenvalue estimates that go
back to work of Weyl \cite{weyl1}, \cite{weyl2}.  Let  $A$ be a 
positive compact operator with eigenvalue list 
$\lambda_1\geq\lambda_2\geq\cdots$ and let 
$\mathcal P_n$ be the set of all $n$-dimensional projections 
in $\mathcal B(H)$.  Then we have 
\begin{equation}\label{EqWeyl}
\sup_{P\in\mathcal P_n}\tr AP =\max_{P\in\mathcal P_n}\tr AP =
\lambda_1+\cdots+\lambda_n,   
\end{equation}
the maximum being achieved on any $n$-dimensional projection 
whose range contains eigenvectors for $\lambda_1,\dots,\lambda_n$.  
Ky Fan's version of this result can be found 
on p. 22 of \cite{Bhatia}.  

The proof of Theorem \ref{thm2} also requires a geometric result, 
asserting that if $p=(p_1,\dots,p_n)$ and 
$\lambda=(\lambda_1,\dots,\lambda_n)$ are two 
finite eigenvalue lists 
that satisfy the first $n$ inequalities (\ref{Eq1}), then 
the components of $\lambda$ can be reduced so as to preserve 
the first $n-1$ inequalities, with equality 
in the $n$th.  

\begin{lem}\label{lem1}
Let $p=(p_1,\dots,p_n)$ and $\lambda=(\lambda_1,\dots,\lambda_n)$ 
be two decreasing sequences of nonnegative reals 
of length $n=1,2,\dots$ satisfying 
\begin{equation}\label{Eq3}
p_1+\cdots+p_k\leq \lambda_1+\cdots+\lambda_k, 
\qquad k=1,2,\dots,n.  
\end{equation} 
There is a decreasing sequence $\mu=(\mu_1,\dots,\mu_n)$ 
such that 
\begin{equation}\label{Eq4}
0\leq \mu_k\leq \lambda_k,\qquad p_1+\cdots+p_k\leq \mu_1+\cdots+\mu_k,
\end{equation}
for $1\leq k\leq n$, 
and $p_1+\cdots+p_n=\mu_1+\cdots+\mu_n$.  
\end{lem}

\begin{proof}
We argue by induction, the case $n=1$ being obvious.   
Fix $n\geq 2$ and suppose that Lemma \ref{lem1} is true for sequences 
of length $n-1$.  Let $D$ be the set of all 
points $\mu=(\mu_1,\dots,\mu_n)\in \mathbb R^n$ satisfying 
$\mu_1\geq\cdots\geq \mu_n\geq 0$ and 
$\mu_k\leq \lambda_k$, $1\leq k\leq n$, and consider 
the compact convex set $K\subseteq \mathbb R^n$ 
$$
K=\{\mu\in D:  \mu_1+\cdots+\mu_k\geq p_1+\cdots+p_k,
\quad k=1,\dots,n-1\}.  
$$
Since 
$f(x)=x_1+\dots+x_n$ is a linear functional on $\mathbb R^n$, 
$f(K)$ is a closed interval $I\subseteq \mathbb R$.  We have to show  
that $p_1+\cdots+p_n\in I$.  For that, it suffices 
to show that there are points $x,y\in K$ such that 
$f(x)\leq p_1+\cdots+p_n\leq f(y)$.  Setting $y=\lambda\in K$, 
we have $p_1+\cdots+p_n\leq \lambda_1+\cdots+\lambda_n=f(y)$ 
by (\ref{Eq3}).  For $x$, use the induction hypothesis 
to obtain numbers $\mu_1\geq\cdots\geq \mu_{n-1}\geq 0$ satisfying 
$0\leq \mu_k\leq \lambda_k$, 
$\mu_1+\cdots+\mu_k\geq p_1+\cdots+p_k$, $1\leq k\leq n-1$, and 
$\mu_1+\dots+\mu_{n-1}=p_1+\cdots+p_{n-1}$.  The point 
$x=(\mu_1,\dots,\mu_{n-1},0)$ belongs to $K$ and 
satisfies 
$
f(x)=p_1+\cdots+p_{n-1}\leq p_1+\cdots+p_n.
$
\end{proof}

\begin{proof}[Proof of Theorem \ref{thm2}]

(i)$\implies$(ii):  Let $e_1, e_2,\dots$ be 
an orthonormal basis for 
$H$ with the property that $\langle Be_j,e_j\rangle=p_k$, 
$j=1,2,\dots$.  Fixing $k$ and letting $E$ be the projection onto 
the span of $e_1,\dots,e_k$, we have 
$$
p_1+\cdots+p_k=\tr(BE)=\tr(L^*ALE)=\tr(ALEL^*)\leq \tr(AF)
$$
where $F$ is the projection onto the range of 
the positive contraction 
$LEL^*$.  Since $F$ is a projection of rank at most 
$k$, the estimate (\ref{EqWeyl}) implies 
$$
\tr(AF)\leq \sup_{\dim F=k}\tr(AF)=\lambda_1+\cdots+\lambda_k, 
$$
and (ii) follows.  

(ii)$\implies$(i):  Let $B$ be a positive compact 
operator in $\mathcal A$ whose eigenvalue 
list $p=(p_1\geq p_2\geq\cdots)$ 
satisfies the inequalities (ii) and let 
$e_1, e_2,\dots$ be an orthonormal basis for $H$ 
such that $[e_1], [e_2],\dots$ are the minimal projections 
of $\mathcal A$. 
Since every permutation of the basis $\{e_k\}$ is 
implemented by a unitary operator $W\in\mathcal B(H)$ satisfying 
$W\mathcal A W^*=\mathcal A$, we may assume without 
essential loss that $ B e_k =p_ke_k$, 
$k=1,2,\dots$.  

We 
construct a sequence of operators $L_n\in\mathcal B(H)$, 
$n=1,2,\dots$,  as follows.   Consider the spectral 
representation of $A$ 
$$
A=\sum_{k=1}^\infty \lambda_k \,\xi_k\otimes\bar\xi_k
$$
where $\xi_1,\xi_2,\dots$ is an orthonormal sequence in 
$H$ consisting of eigenvectors of $A$.  Fix $n$, let 
$H_n$ be the linear span of $\xi_1,\dots,\xi_n$, and 
let $A_n$ be the restriction of $A$ to $H_n$.  The eigenvalue list  
of $A_n$ is $(\lambda_1,\dots,\lambda_n)$; so by Lemma \ref{lem1}, 
there is a decreasing sequence $\mu=(\mu_1,\dots,\mu_n)$ satisfying 
$0\leq \mu_k\leq \lambda_k$ for $1\leq k\leq n$, and 
$$
p_1+\cdots+p_k\leq \mu_1+\cdots+\mu_k, \qquad k=1,\dots, n,
$$  
with equality holding for $k=n$.  
The sequence $\mu$ of course 
depends on $n$ but we suppress that 
in the notation since $\mu$ will soon disappear.  Consider 
the operator $B_n$ defined on $H_n$ by requiring 
$B_n\xi_k=\mu_k\xi_k$, $1\leq k\leq n$.  The eigenvalue list of 
$B_n$ dominates $(p_1,\dots,p_n)$ as 
in the hypothesis of Horn's result Theorem \ref{hornThm}.  
Thus there is an orthonormal basis $e_1^{(n)}, \dots,e_n^{(n)}$ for 
$H_n$ such that 
\begin{equation*}
\langle B_ne_k^{(n)},e_k^{(n)}\rangle=p_k,\qquad 
k=1,\dots,n.  
\end{equation*}
Since $0\leq B_n\leq A_n$ it follows that 
$$
p_k\leq \langle A_ne_k^{(n)},e_k^{(n)}\rangle=
\langle Ae_k^{(n)},e_k^{(n)}\rangle, 
\qquad k=1,\dots,n.  
$$
Let $L_n\in\mathcal B(H)$ be the operator defined by 
$L_ne_k= e_k^{(n)}$ for $k=1,\dots,n$, and $L_n=0$ on 
the orthocomplement of $[e_1^{(n)},\dots,e_n^{(n)}]$.  

We have constructed a sequence $L_1, L_2,\dots$ of finite rank 
partial isometries in $\mathcal B(H)$ 
that satisfies the system of inequalities 
\begin{equation}\label{Eq6}
p_k\leq \langle AL_ne_k,L_ne_k\rangle,\qquad n\geq k\geq 1.  
\end{equation}
Since the unit ball of $\mathcal B(H)$ is sequentially 
compact in its weak operator topology, there is a 
subsequence $n_1<n_2<\dots$ 
and a contraction $L_\infty\in\mathcal B(H)$ such that 
$\langle L_{n_j}\eta,\zeta\rangle\to\langle L_\infty \eta,\zeta\rangle$ 
as $j\to\infty$,
for every $\eta,\zeta\in H$.  
We claim that $L_\infty$ 
satisfies 
\begin{equation}\label{Eq7}
p_k\leq \langle AL_\infty e_k,L_\infty e_k\rangle, \qquad k=1,2,\dots.  
\end{equation}
To see that, fix $k$ and note that for sufficiently 
large $j$, (\ref{Eq6}) implies 
$$
p_k\leq \langle AL_{n_j}e_k,L_{n_j}e_k\rangle.  
$$
As $j\to\infty$, $L_{n_j}e_k$ tends to $L_\infty e_k$ in the weak 
topology of $H$.  Since $A$ is a compact operator, 
$\|AL_{n_j}e_k- AL_\infty e_k\|\to 0$ as $j\to\infty$; hence the 
inner products $\langle AL_{n_j}e_k,L_{n_j}e_k\rangle$ 
converge to $\langle AL_\infty e_k,L_\infty e_k\rangle$, 
and (\ref{Eq7}) follows.  

Finally, choose $t_1, t_2,\dots\in[0,1]$ such that 
$p_k=t_k\langle A L_\infty u_k, L_\infty u_k\rangle$ 
for every $k$.  
Letting $D\in\mathcal B(H)$ be the contraction defined by 
$De_k=\sqrt{t_k}e_k$, $k\geq 1$, one finds that 
the operator $L=L_\infty D$ satisfies 
\begin{equation*}
\langle L^*A  L e_k,e_k\rangle = p_k,
\qquad k=1,2,\dots,  
\end{equation*}
and the required formula 
$$
E(L^*AL)=\sum_{k=1}^\infty E_k L^*AL E_k=\sum_{k=1}^\infty p_k E_k=B
$$
follows.  
\end{proof}

\begin{proof}[Proof of Theorem \ref{thm1}]  Let $E_1, E_2,\dots$ 
be the minimal projections of $\mathcal A$ and let 
$e_1, e_2,\dots$ be an orthonormal basis for $H$ such 
that $E_k$ is the projection $[e_k]$, 
$k=1,2,\dots$.  

We show first that for every positive trace class operator 
$A\in\mathcal B(H)$ with eigenvalue list $\lambda$, 
the eigenvalue list $p=(p_1,p_2,\dots)$ of 
$B=E(A)$ must satisfy (\ref{Eq1}) and (\ref{Eq2}).  By 
permuting the elements of the 
basis $\{e_k\}$ appropriately and changing notation, 
we may assume that $Be_k=p_ke_k$, $k=1,2,\dots$.  
Let $P_n$ be the projection on $[e_1,\dots,e_n]$.  
Since $A$ is a 
positive compact operator with eigenvalue list $\lambda$, 
we can make use of (\ref{EqWeyl}) to write 
$$
p_1+\cdots+p_n=\sum_{k=1}^n\langle Be_k,e_k\rangle=
\tr AP_n\leq \lambda_1+\cdots+\lambda_n.     
$$
Moreover, 
$p_1+p_2+\cdots = \langle Au_1,u_1\rangle+\langle Au_2,u_2\rangle+\cdots
= \tr A=\lambda_1+\lambda_2+\cdots$.   

Conversely,  let $p$ and $\lambda$ be two summable eigenvalue lists that satisfy 
(\ref{Eq1}) -- (\ref{Eq2}), and let $B$ be a positive  
trace-class operator in $\mathcal A$ with list $p$.  Again, by relabeling 
the orthonormal basis $\{e_k\}$, we may assume that 
$Be_k=p_ke_k$, $k=1,2,\dots$.  Choose any positive trace-class operator 
$A\in\mathcal B(H)$ having eigenvalue list $\lambda$, and let 
$P$ be the projection onto the closure of $AH$.  
Theorem \ref{thm2} implies that there is a contraction 
$L\in\mathcal B(H)$ satisfying 
$p_k= \langle ALe_k,Le_k\rangle$ for $k\geq 1$.  
By replacing $L$ with $PL$ if necessary, we may 
also assume that $LH$ is contained in $PH$, and in 
that case we claim: 
\begin{equation}\label{Eq9}
LL^* =  P.  
\end{equation}
Indeed, 
$P-LL^*\geq 0$ because $L$ is a contraction whose 
range is contained in $PH$, 
and it suffices to show that 
the positive operator $A^{1/2}(P-LL^*)A^{1/2}=A-A^{1/2}LL^* A^{1/2}$ 
has trace zero; equivalently, 
$\tr A^{1/2}LL^*A^{1/2} = \tr A$.   Using
$\tr XX^*=\tr X^*X$ for $X=A^{1/2}L$, we have 
\begin{align*}
\tr A^{1/2}LL^*A^{1/2}&=\tr L^*AL=\sum_{n=1}^\infty \langle ALe_n,Le_n\rangle 
\\
&=p_1+p_2+\cdots
=\lambda_1+\lambda_2+\cdots=\tr A.   
\end{align*}

$L$ is a co-isometry by (\ref{Eq9}); hence it can be 
changed into a unitary operator $U: H\to H\oplus \ker L$ by 
making use of the projection $Q: H\to \ker L$ as follows: 
$U\xi=L\xi\oplus Q\xi$, $\xi\in H$.  
Now consider the operator 
$$
A_0\oplus \mathbf 0\in \mathcal B(PH\oplus\ker L),    
$$  
$A_0$ denoting the restriction of $A$ to $PH=\overline{AH}$.  
Since 
$Ue_k=Le_k\oplus Q e_k$, $k=1,2,\dots$, we have 
\begin{equation*}
\langle (A_0\oplus\mathbf 0) Ue_k,Ue_k\rangle =\langle ALe_k,Le_k\rangle =p_k,
\qquad k=1,2,\dots.    
\end{equation*}
Therefore $U^*(A_0\oplus\mathbf 0)U$ is a 
positive trace class operator in $\mathcal B(H)$ satisfying 
$$
E(U^*(A_0\oplus\mathbf 0)U)=\sum_{k=1}^\infty p_k E_k=B.  
$$
Since $U^*(A_0\oplus\mathbf 0)U$ has the same eigenvalue 
list as $A$, Proposition \ref{prop1} implies that it must 
belong to $\mathcal O(A)=\mathcal O_\lambda$, and the proof is complete.  
\end{proof}


In the series \cite{kadPnasI}, \cite{kadPnasII}, one of us  
carried out a study of the possible diagonals of projections acting 
on a separable Hilbert space.  
The results of this paper 
do not address the most difficult case where the projection 
has infinite rank and infinite co-rank; but Theorem \ref{thm1} 
does give  
the result of Theorem 13 of \cite{kadPnasII}, as follows.  We formulate 
that in terms of the Hilbert space $\ell^2$ and its standard 
orthonormal basis $(u_k)_{k\geq1}$, with the associated  
realization of operators as matrices relative to this basis.

\begin{cor}\label{kadThm}
Let  $p=(p_1,p_2,\dots)$ be a sequence of numbers in the unit interval 
$0\leq p_k\leq 1$, and let $m$ be a positive integer.  The following 
are equivalent
\begin{enumerate}
\item[(i)]
There is a rank $m$ projection $P\in\mathcal B(\ell^2)$ whose matrix 
has $p$ as its diagonal.  
\item[(ii)]
$p_1+p_2+\cdots=m$.  
\end{enumerate}  
\end{cor}

\begin{proof} We prove the nontrivial implication 
(ii)$\implies$ (i).  Since (ii) implies that the sequence 
$p_n$ converges to zero and since permutations of $\mathbb N$ are 
implemented by unitary operators on $\ell^2(\mathbb N)$ in the obvious 
way, it suffices to address the case where the sequence 
is decreasing $p_1\geq p_2\geq\cdots$.  The eigenvalue 
list of a projection of rank $m$ is 
$$
\lambda=(\underbrace{1,\dots,1}_{m{\text\ times}},0,0,\dots), 
$$
and $\mathcal O_\lambda$ consists of all rank $m$ projections
in $\mathcal B(\ell^2)$.  
The hypothesis (ii), together with $0\leq p_k\leq 1$, implies 
that 
$
p_1+\cdots+p_n\leq\lambda_1+\cdots+\lambda_n
$
holds for every 
$n\geq 1$. Hence Theorem 
\ref{thm1} implies that there is an operator 
in $\mathcal O_\lambda$ with diagonal sequence $p$.  
\end{proof}

\part{Type $II_1$ Factors}

We turn now to the case of self-adjoint operators $A$ in 
a finite factor.  In this context, the appropriate 
counterpart of the eigenvalue list is a finite positive 
measure on the real line, called the spectral distribution 
of $A$.  After working out the basic properties of 
the spectral distribution and relating it to 
approximate unitary equivalence, we establish 
a generalization of the Schur inequalities for 
$II_1$ factors.

\section{Spectral distribution of a self-adjoint operator}\label{S:specDist}
In the remainder of this paper we work within the context 
of a {\em finite} factor 
$\mathcal M$ with 
normalized trace $\tau$.  For every self-adjoint operator $A\in\mathcal M$ 
there is a unique  probability measure $m_A$ 
on the Borel subsets of the real line whose moments are given by 
\begin{equation}\label{Eq11}
\int_{-\infty}^\infty \lambda^n\,dm_A(\lambda)=\tau(A^n),\qquad n=0,1,2,\dots.  
\end{equation}
\begin{defn}
The measure $m_A$ is called the spectral distribution of $A$.  
\end{defn}
The purpose of these section is to discuss the basic  
properties of this invariant.  
The spectral distribution is the appropriate generalization 
to $II_1$ factors of the 
eigenvalue list invariant
of self-adjoint $n\times n$ matrices.  Indeed, if 
$A$ is a self-adjoint $n\times n$ matrix with eigenvalue list 
$\Lambda=\{\lambda_1\geq\lambda_2\geq\cdots\}$, 
then  $m_A$ is the discrete measure 
$$
m_A=\frac{1}{n}(\delta_{\lambda_1}+\cdots+\delta_{\lambda_n}),   
$$
$\delta_\lambda$ denoting the unit point mass 
concentrated at $\lambda\in\mathbb R$.  
Equivalently, $m_A$ 
assigns mass to singletons $\{\lambda\}$ of $\mathbb R$ as follows
$$
m_A({\lambda})=
\begin{cases}
\frac{1}{n}({\rm multiplicity\,of }\ \lambda),\quad& {\rm if\,} \lambda\in \sigma(A)\\
0,& {\rm otherwise}.  
\end{cases}
$$
We require the following observation, 
which asserts that the spectral distribution 
of an operator in a $II_1$ factor can be arbitrary.

\begin{prop}\label{prop0}
Let $\mathcal A\subseteq\mathcal M$ be a MASA in a 
$II_1$ factor and let $m$ be a compactly supported probability 
measure on the real line.  Then there is a self-adjoint operator 
$A\in\mathcal A$ such that $m=m_A$.  
\end{prop}  

\begin{proof}
$\mathcal A$ contains a countably-generated nonatomic subalgebra, 
which must be isomorphic to $L^\infty[0,1]$ in such a way that 
the restriction of the trace $\tau$ corresponds to the state 
of $L^\infty[0,1]$ given by 
$$
\tau(f)=\int_0^1 f(x)\,dx, \qquad f\in L^\infty[0,1].  
$$
Thus it suffices to show that there is a real-valued function 
$f\in L^\infty[0,1]$ such that $\int_0^1 f(x)^n\,dx=\int \lambda^n\,dm(\lambda)$ 
for $n=0,1,\dots$ or equivalently, for every Borel set $S\subseteq\mathbb R$, 
\begin{equation}
\tau\{x\in[0,1]:f(x)\in S\} = m(S),   
\end{equation}
where we abuse notation slightly by also writing $\tau$ for Lebesgue 
measure on the unit interval.  Let $K$ be the closed 
support of $m$.  The pair $(K,m)$ defines a separable measure algebra 
which may have a finite or countable number of atoms.  On the other 
hand, $([0,1],\tau)$ gives rise to a nonatomic separable measure algebra.  
Letting $x_1, x_2,\dots$ be the points of $K$ having positive 
$m$-measure, we find a sequence of disjoint Borel sets 
$E_1, E_2,\dots\subseteq[0,1]$ such that $\tau(E_k)=m(\{x_k\})$ 
for all $k$.  Define $f$ on $\cup_k E_k$ so that it takes 
on the constant value $x_k$ throughout $E_k$.  $f$ is 
a measure-preserving map of $\cup_k E_k$ onto the atomic 
part of $(K,\tau)$.  Since 
$\tau([0,1]\setminus\cup_kE_k)=m(K\setminus \{x_1,x_2,\dots\})$ 
and since the remaining parts of both measure spaces are nonatomic 
and separable, there is a surjective Borel map of 
$[0,1]\setminus \cup_n E_n$ onto $K\setminus \{x_1,x_2,\dots\}$ 
that pushes Lebesgue measure forward to $m$, and we can use this 
map to extend the definition of $f$ to all of $[0,1]$ in the obvious way.  
\end{proof}

The eigenvalue list is a complete invariant 
for $\mathcal L^1$-equivalence of positive trace-class operators in 
type $I$ factors.   
We now show that the spectral distribution 
invariant occupies a similar position.  
\begin{defn}\label{Def1}
Operators $A,B\in\mathcal M$ are said to 
be {\em approximately equivalent} if there is a sequence 
of unitary operators $U_1, U_2,\dots$ in  $\mathcal M$ 
such that 
\begin{equation}\label{Eq33}
\lim_{n\to\infty}\|U_nAU_n^*-B\|=0.  
\end{equation}
The set of all operators  in the norm-closed unitary 
orbit of $A$ is written $\mathcal O_A$.  
\end{defn}
\begin{thm}\label{thm11}
Let $A$ be a self-adjoint operator in $\mathcal M$,  let 
$m_A$ be the spectral distribution of $A$, and let 
$\mathcal O_A$ be the norm-closed unitary orbit of $A$.  
Then $\mathcal O_A$ 
is closed in the strong operator topology, and consists of 
all self-adjoint operators $B\in \mathcal M$ satisfying 
$m_B=m_A$.  
\end{thm}

Before giving the proof we collect an elementary observation.  

\begin{lem}\label{lem11}
Let $ E_1\leq E_2\leq\cdots \leq E_n$ and 
$F_1\leq F_2\leq\cdots \leq F_n$ be two linearly ordered sets of 
projections in $\mathcal M$ such that $\tau(E_k)=\tau(F_k)$ 
for $k=1,\dots,n$.  Then there is a unitary operator 
$U$ in $\mathcal M$ such that $UE_kU^*=F_k$, $k=1,\dots,n$.  
\end{lem}

\begin{proof} By adjoining the identity to the end of 
each list if necessary, 
we can assume that $E_n=F_n=\mathbf 1$.  
Setting $E_0=F_0=0$, the hypothesis implies that 
$\tau(E_k-E_{k-1})=\tau(F_k-F_{k-1})$ for each $k=1,\dots,n$.  Since 
$\mathcal M$ is a finite factor, projections with the same trace 
must be Murray-von Neumann equivalent.  Thus there are partial 
isometries $U_1,\dots, U_n\in \mathcal M$ with 
$U_k^*U_k=E_k-E_{k-1}$ and $U_kU_k^*=F_k-F_{k-1}$ for all $k$.  
The projections $U_1^*U_1,\dots,U_n^*U_n$ add up to $E_n=\mathbf 1$, 
and similarly 
$U_1U_1^*+\cdots+U_nU_n^*=\mathbf 1$.  It follows that 
$W=U_1+\cdots+U_n$ is a unitary operator in $\mathcal M$ with 
the property $W(E_k-E_{k-1})W^*=F_k-F_{k-1}$ for every $k=1,\dots,n$, 
hence $WE_kW^*=F_k$ for $k=1,\dots,n$.  
\end{proof}

\begin{proof}[Proof of Theorem \ref{thm11}]  
We will show that a self-adjoint operator $B$ belongs to 
$\mathcal O_A$ iff $m_B=m_A$.  Once that is established, it 
will follow that $\mathcal O_A$ is strongly closed, because  
the relation $m_B=m_A$ is characterized by the sequence of 
equations 
\begin{equation}\label{Eq3.1}
\tau(B^n)=\tau(A^n),\qquad n=0,1,2,\dots  
\end{equation}
and each of the monomials $B\mapsto \tau(B^n)$ is strongly 
continuous on bounded subsets of $\mathcal M$.  

Every operator $B\in \mathcal M$ that is unitarily 
equivalent to $A$ must satisfy the formulas (\ref{Eq3.1}), 
and hence a norm-limit of such operators will satisfy 
the same formulas.  It follows that $m_B=m_A$ for 
every $B$ in the norm-closed unitary orbit of $A$.  

Conversely, let $A$, $B$ be self-adjoint operators 
of $\mathcal M$ satisfying $m_A=m_B$, and let 
$[a,b]$ be an interval with the property that 
$m_A=m_B$ is supported 
in the interior $(a,b)$.  
Then the spectra of $A$ and $B$ are both contained 
in $(a,b)$.  
For every $a\leq t\leq b$ let $E_t$ (resp. $F_t$) be the spectral projection 
of $A$ (resp. $B$) corresponding to the subinterval 
$[a,t]\subseteq \mathbb R$.  Then by hypothesis we have 
\begin{equation}\label{Eq44}
\tau(E_t)=m_A([a,t])=m_B([a,t])=\tau(F_t), \qquad a\leq t\leq b.  
\end{equation}
Given $\epsilon>0$ we can find a partition 
$a=t_0<t_1<\cdots<t_n=b$ of $[a,b]$ fine 
enough that 
$$
|t-\sum_{k-1}^n t_k\chi_{(t_{k-1},t_k]}(t)|\leq \epsilon,\qquad 
a\leq t\leq b.     
$$
Letting $A_0$, $B_0$ be the operators 
$$
A_0=\sum_{k=1}^nt_k(E_{t_k}-E_{t_{k-1}}),\quad 
B_0=\sum_{k=1}^nt_k(F_{t_k}-F_{t_{k-1}}),   
$$
we find that $\|A-A_0\|\leq \epsilon$ and $\|B-B_0\|\leq \epsilon$.  
(\ref{Eq44}) implies $\tau(E_{t_k})=\tau(F_{t_k})$ for every $k$, so by 
the Lemma there is 
a unitary operator $W\in\mathcal M$ such that 
$WE_{t_k}W^*=F_{t_k}$ for all $k$, hence $WA_0W^*=B_0$.  
An obvious estimate now implies 
$\|WAW^*-B\|\leq 2\epsilon$, and since $\epsilon$ is 
arbitrary it follows that $A$ and $B$ are approximately 
equivalent.   
\end{proof}

Theorem \ref{thm11} implies that the spectral distribution 
is a complete invariant for approximate unitary equivalence, and 
it is natural to ask if two self-adjoint operators 
that are approximately equivalent must be unitarily 
equivalent, or at least conjugate by way of a 
$*$-automorphism.   The following class of 
examples shows that the answer is no.  
\vskip0.1in
{\bf Example.}
Let $\mathcal R$ be a $II_1$ factor and 
let $\mathcal A$ and $\mathcal B$ be two MASAs in 
$\mathcal R$ that are {\em not} conjugate by way of an 
automorphism of $\mathcal R$.  For example, $\mathcal A$ 
can be taken to be a regular MASA and $\mathcal B$ a 
singular one.  Since both $\mathcal A$ and $\mathcal B$ 
are isomorphic to $L^\infty[0,1]$ by way of an isomorphism 
that carries the trace to Lebesgue measure, it follows 
that a) there is a $*$-isomorphism $\alpha$ of $\mathcal A$ 
onto $\mathcal B$  satisfying $\tau(\alpha(X))=\tau(X)$ 
for all $X\in\mathcal A$, and b) $\mathcal A$ is the 
von Neumann algebra $W^*(A)$ generated by a single self-adjoint 
operator $A$.  Let $B=\alpha(A)$.  Then $\tau(A^n)=\tau(B^n)$ 
for every $n=0,1,2,\dots$ and hence $m_A=m_B$.  It follows 
from Theorem \ref{thm11} that 
$A$ and $B$ are approximately equivalent.  On the 
other hand, there is no unitary operator $U\in\mathcal R$ 
satisfying $UAU^*=B$, since that would imply that 
$\theta(X)=UXU^*$ is an automorphism of $\mathcal R$ 
that carries $\mathcal A=W^*(A)$ onto $\mathcal B=W^*(B)$.  
\vskip0.1in

\section{Schur-type inequalities for $II_1$ factors}\label{S:sec5}

The purpose of this section is to formulate an 
appropriate counterpart of 
the Schur inequalities for self-adjoint 
operators in $II_1$ factors. This is not 
the only formulation possible, and we refer the reader 
to Section 5 of \cite{kadSchur} for an alternate approach.  
Here, we seek to formulate the Schur inequalities in terms 
of spectral distributions.  
That formulation is based on the 
following observations.  

\begin{prop}\label{propII1}
For any two compactly supported probability measures $m$, $n$
on the real line, 
the following are equivalent:
\begin{enumerate}
\item[(i)]$m$ and $n$ have the same first moment
$$
\int_\mathbb R \lambda\,dm(\lambda)=\int_\mathbb R \lambda\,dn(\lambda),  
$$
and for every $t\in\mathbb R$ we have 
$$
\int_t^\infty m([s,\infty))\,ds\leq \int_t^\infty n([s,\infty))\,ds.  
$$
\item[(ii)]
$m$ and $n$ have the same first moment,
and for every $t\in\mathbb R$ we have  
$$
\int_{[t,\infty)}(\lambda-t)\,dm(\lambda)\leq 
\int_{[t,\infty)}(\lambda-t)\,dn(\lambda).   
$$
\item[(iii)]For every continuous convex function defined on 
a closed 
interval $I=[c,d]$ that supports both measures $m$ and $n$, we have 
$$
\int_I f(\lambda)\,dm(\lambda)\leq \int_I f(\lambda)\,dn(\lambda).  
$$
\end{enumerate}
\end{prop}

\begin{proof}
The equivalence (i)$\iff$(ii) is a consequence of the 
classic integration by parts formula 
of Riemann-Stieltjes integration, which can be applied as follows.    
Fix $t\in \mathbb R$, 
let $m$ be a compactly supported measure defined on $\mathbb R$, 
and choose $a,b\in\mathbb R$ so that 
$a<t<b$ and such that $(a,b)$ contains the closed 
support of both $m$ and $n$.  
Let $\alpha:\mathbb R\to\mathbb R$ be 
the decreasing function $\alpha(s)=m([s,\infty))$ 
and let $f$ be the continuous increasing function $f(s)=\max(s-t,0)$.  
An application of Theorem 9--6 of \cite{apostolBook} gives 
$$
\int_a^bf(x)\,d\alpha(x)+\int_a^b\alpha(x)\,df(x)=f(b)\alpha(b)-f(a)\alpha(a).  
$$
In this case, 
$\alpha(b)=f(a)=0$, and straightforward computations show that 
\begin{align*}
\int_a^bf(x)\,d\alpha(x)&=-\int_{[t,\infty]}(x-t)\,dm(x), 
\\
\int_a^b\alpha(x)\,df(x)&=\int_t^\infty m([x,\infty))\,dx.  
\end{align*}
It follows that 
$$
\int_{[t,\infty]}(\lambda-t)\,dm(\lambda)=\int_t^\infty m([s,\infty))\,ds, 
$$
and the equivalence of (i) and (ii) follows.

(ii)$\implies$(iii):  This will follow 
if we show that every 
continuous real-valued convex function $f$ defined on a compact 
interval $I\subseteq \mathbb R$ 
can be uniformly approximated on 
$I$ by functions 
\begin{equation}\label{Eq01}
f(\lambda)=a+b\lambda+g(\lambda)
\end{equation}
where $a$ and $b$ are real constants and 
$g$ belongs to the cone generated by the ``angular" functions 
$$
g_t(\lambda)=\max(\lambda-t,0)=
(\lambda-t)\chi_{[t,\infty)}(\lambda),\qquad t\in\mathbb R.
$$ 
To see how the approximation (\ref{Eq01}) is 
achieved, one first 
approximates $f$ uniformly on $I$ 
with a twice continuously 
differentiable convex function $g$.  
Since $g^\prime$ is an increasing function, it can be 
uniformly approximated by an increasing step 
function having the form $a + h(\lambda)$ 
where $a$ is constant and $h$ belongs to the cone generated by 
the step functions $\chi_{[t,\infty)}$, $t\in\mathbb R$.  After 
one integrates this approximation of $g^\prime$ 
one obtains an approximation to $g(\lambda)$ of the form 
$a\lambda+b+\int h(\lambda)\,d\lambda$.  
Moreover, since the indefinite integral of a step 
function $\chi_{[t,\infty)}(\lambda)$ 
has the form $c+g_t(\lambda)$ where $c$ is 
a constant, this 
approximation of $g$ has the 
form (\ref{Eq01}).  

The implication (iii)$\implies$(ii) is an immediate 
consequence of the fact 
that functions of the form (\ref{Eq01}) are continuous 
and convex.  
\end{proof}

\begin{defn}
Let $m$ and $n$ be two compactly supported probability measures on 
the real line.  We write $m\preceq n$ if $m$ is dominated by $n$ 
in the equivalent senses of Proposition \ref{propII1}, i.e., if 
$\int_\mathbb R\lambda\,dm(\lambda)=\int_\mathbb R\lambda\,dn(\lambda)$, 
and 
\begin{equation}\label{schIneq}
\int_t^\infty m([s,\infty))\,ds\leq 
\int_t^\infty n([t,\infty))\,ds, \qquad t\in\mathbb R.  
\end{equation}
\end{defn}

The relation $\preceq$ is obviously a partial ordering on the 
set of all compactly supported probability measures on the real 
line.  Given two self-adjoint operators $A,B$ in a $II_1$ factor 
$\mathcal M$, we interpret the relation 
$m_A\preceq m_B$ as the appropriate counterpart of the 
Schur inequalities (\ref{Eq000}) that relate the eigenvalue 
lists of $A$ and $B$.  
This interpretation is justified by Proposition 
\ref{propII1} and the following remarks.  

\begin{rem}[Relation to the classical inequalities of Schur]
Let $\tau$ be the normalized trace on the matrix algebra 
$M_n(\mathbb C)$, and let $A$ and $B$ be  self-adjoint 
$n\times n$ matrices.  We have discussed the relation between 
the eigenvalue list of $A$ and the spectral distribution 
$m_A$ in Section \ref{S:specDist}.  We now 
examine the relation between the integrals
$$
\int_t^\infty m_A([s,\infty))\,ds 
$$ 
appearing in 
(\ref{schIneq}) and the eigenvalue list 
$\lambda_1\geq\cdots\geq \lambda_n$ 
of $A$.  For simplicity, we consider the 
case where the eigenvalues are simple ones.  For every $t$ in the 
interval $\lambda_{k+1}<t<\lambda_k$ one has $m([t,\infty))=k/n$, 
and after a straightforward integration and cancellation one 
finds that for $\lambda_{k-1}<t<\lambda_k$,  
$$
\int_t^\infty m_A([s,\infty))\,ds = \frac{\lambda_1+\lambda_2+\cdots+\lambda_k-kt}{n}.  
$$
Let $B$ be another self-adjoint matrix with eigenvalue 
list $\mu_1\geq\cdots\geq \mu_n$.  
The preceding formula shows that the system of inequalities 
\begin{equation}\label{SchurEq2}
\int_t^\infty m_A([s,\infty))\,ds\leq \int_t^\infty m_B([s,\infty))\,ds, 
\qquad t\in\mathbb R, 
\end{equation}
differs somewhat from the system of 
classical Schur inequalities, which in terms of the normalized trace 
would assert 
\begin{equation}\label{SchurEq33}
\frac{\lambda_1+\cdots+\lambda_k}{n}\leq \frac{\mu_1+\cdots+\mu_k}{n},\qquad k=1,2,\dots,n.  
\end{equation}
However, if $\tau(A)=\tau(B)$ then 
$\lambda_1+\cdots+\lambda_n=\mu_1+\cdots+\mu_n$; and in 
that event  
the inequalities (\ref{SchurEq2}) are {\em equivalent} to the 
Schur inequalities (\ref{SchurEq33}) because they are equivalent to the 
inequalities of assertion (iii) of Proposition \ref{propII1}.  
That is a consequence of classical results of Hardy, Littlewood and 
Polya \cite{hlp} which are summarized in Theorem 1 of \cite{hornDiag}.  
The relevant result asserts that 
for two finite eigenvalue lists 
$$
\{\lambda_1\geq\cdots\geq\lambda_n\},\quad \{\mu_1\geq\cdots\geq\mu_n\}
$$
which satisfy $\lambda_1+\cdots+\lambda_n=\mu_1+\cdots+\mu_n$, 
the following are equivalent:
\begin{enumerate}
\item[(1)]$\lambda_1+\cdots+\lambda_k\leq \mu_1+\cdots+\mu_k$, for every 
$k=1,\dots,n$.  
\item[(2)]For every convex function $f$ defined on an interval containing 
all $\lambda_i$ and $\mu_j$, one has 
$$
\sum_{k=1}^nf(\lambda_k)\leq \sum_{k=1}^n f(\mu_k).  
$$
\end{enumerate} 
Thus, when taken together with the equivalence of (1) 
and (2), Proposition \ref{propII1} implies that 
{\em the system of inequalities (\ref{schIneq}) is an appropriate generalization 
of the Schur inequalities (\ref{SchurEq33}) to $II_1$ factors. } 
\end{rem}

\section{Proof of the Schur inequalities}\label{S:sec6}

We require a convexity inequality for operators 
in a $II_1$ factor.  While related results can be 
found in the literature, we have been unable to find 
references appropriate 
for this particular result, and we include 
a proof for completeness.  Let $\mathcal A$ be a maximal abelian self-adjoint 
subalgebra of a 
$II_1$ factor $\mathcal M$ with normalized trace $\tau$, and 
let $E: \mathcal M\to \mathcal A$ be the $\tau$-preserving 
conditional expectation.  

\begin{prop}\label{prop1.1}
Let $f$ be a real-valued continuous convex function 
defined on a compact interval $I=[a,b]\subseteq \mathbb R$.  
Then for every self-adjoint operator $A\in\mathcal M$ with 
spectrum contained in $I$, the spectrum of $E(A)$ is also 
contained in $I$ and  we have 
\begin{equation}\label{Eq5}
f(E(A))\leq E(f(A)).  
\end{equation}
\end{prop}

\begin{proof}
For every self-adjoint operator $A\in\mathcal M$, we 
write $A_+$ for the positive part of $A$, defined by 
$A_+=AP_+=P_+A$ where $P_+$ is the spectral projection 
of $A$ associated with the nonnegative real axis 
$[0,\infty)$.  
We claim first that
\begin{equation}\label{Eq5.1}
E(A)_+\leq E(A_+).  
\end{equation}
Indeed, in the natural ordering of 
self-adjoint operators in $\mathcal M$ we have 
$A\leq A_+$ and hence $E(A)\leq E(A_+)$.  Thus $E(A_+)$ is 
a positive operator that dominates $E(A)$.  Since $\mathcal A$ 
is abelian, $E(A)_+$ is the smallest 
positive operator in $\mathcal A$ 
that dominates $E(A)$, and (\ref{Eq5.1}) follows.

In order to prove (\ref{Eq5}), choose $r\in\mathbb R$ and 
let $g_r(\lambda)=\max(\lambda-r,0)$.  We may 
apply (\ref{Eq5.1}) to the operator $A-r\mathbf 1$ to obtain 
$$
g_r(E(A))=(E(A)-r\mathbf 1)_+=E(A-r\mathbf 1)_+\leq E((A-r\mathbf 1)_+)=E(g_r(A)).  
$$
It follows that for every convex 
function $f_0:\mathbb R\to \mathbb R$ 
of the form 
\begin{equation}\label{Eq5.2}
f_0(\lambda)=a+b\lambda+\sum_{k=1}^n c_k g_{r_k}(\lambda)
\end{equation}
where $a,b, r_1,\dots,r_n\in\mathbb R$ and $c_1,\dots,c_n\geq0$, one has 
$$
f_0(E(A))\leq E(f_0(A)).  
$$
Since every continuous convex function $f:[a,b]\to\mathbb R$ 
can be uniformly approximated by functions $f_0$ of 
the form (\ref{Eq5.2}), one deduces 
(\ref{Eq5}) for continuous convex functions from 
these inequalities.  \end{proof}

The following result establishes the 
Schur inequalities for operators in a $II_1$ factor.  
\begin{thm}\label{thm22}
Let $\mathcal A$ be 
a MASA in $\mathcal M$ and let 
$E: \mathcal M\to \mathcal A$ be the trace-preserving 
conditional expectation.  For every self-adjoint 
operator $A$ in $\mathcal M$, the spectral 
distribution of $B=E(A)$ is related to 
that of $A$ by $m_B\preceq m_A$.  
\end{thm}

\begin{proof}
Let $[a,b]$ be the smallest closed interval containing 
$\sigma(A)\cup\sigma(B)$.  
Since both $m_B$ and $m_A$ are probability measures, 
Proposition \ref{prop1} implies that $m_B\preceq m_A$ iff 
for every continuous convex function $f\in C[a,b]$,  
$$
\int_I f(\lambda)\,dm_B(\lambda)\leq \int_I f(\lambda)\,dm_A(\lambda).  
$$
Since the left side is $\tau(f(B))=\tau(f(E(A)))$ and the right 
side is $\tau(f(A))$, the preceding inequality follows 
from formula (\ref{Eq5}).  
\end{proof}

Theorem \ref{thm22} makes the following assertion 
about the norm-closed unitary orbit $\mathcal O_A$ of a self-adjoint 
operator: $E(\mathcal O_A)$ is contained in 
the set of all self-adjoint operators $B\in\mathcal A$ satisfying 
$m_B\preceq m_A$.  Thus, an affirmative reply to 
the following question would appear to be a 
natural generalization of Horn's 
Theorem for $n\times n$ matrices.   
\vskip0.1in
{\bf Problem.}  Let $\mathcal A$ be 
a MASA in a $II_1$ factor $\mathcal M$, let 
$E: \mathcal M\to \mathcal A$ be the trace-preserving 
conditional expectation, and let $A$ be a self-adjoint 
operator in $\mathcal M$.  Does $E(\mathcal O_A)$ contain  
the set of all self-adjoint operators $B\in\mathcal A$ 
satisfying $m_B\preceq m_A$?   
\vskip0.1in
%


\bibliographystyle{alpha}

\newcommand{\noopsort}[1]{} \newcommand{\printfirst}[2]{#1}
  \newcommand{\singleletter}[1]{#1} \newcommand{\switchargs}[2]{#2#1}

\end{document}